\documentclass[12pt]{amsart}
\usepackage{amsmath,amsfonts,amsthm,amssymb,mathrsfs,amsxtra,amscd,latexsym, xcolor}
\usepackage[hmargin=2cm,vmargin=3cm]{geometry}
\usepackage{hyperref}
\usepackage{graphicx}

\usepackage{amsmath,amsthm,amssymb,enumerate,amsfonts,mathrsfs,amsxtra,amscd,latexsym}

\newtheorem{theorem}{Theorem}[section]

\newtheorem{definition}[theorem]{Definition}
\newtheorem{example}[theorem]{Example}

\newtheorem{question}[theorem]{Question}
\newtheorem{questions}[theorem]{Questions}
\newtheorem*{remark}{Remark}

\numberwithin{equation}{section}

\newcommand{\sm}{\left(\begin{smallmatrix}}
\newcommand{\esm}{\end{smallmatrix}\right)}
\newcommand{\mat}{\left(\begin{matrix}}
\newcommand{\emat}{\end{matrix}\right)}

\def\Gal{\mathop{\rm Gal}}

\def\Q{\mathbb Q}  \def\Z{\mathbb Z}

\newcommand{\F}{\mathbb{F}}

\newcommand{\beq}{\begin{eqnarray*}}
\newcommand{\eeq}{\end{eqnarray*}}
\newcommand{\beqn}{\begin{eqnarray}}
\newcommand{\eeqn}{\end{eqnarray}}

\newcommand{\ben}{\begin{enumerate}}
\newcommand{\een}{\end{enumerate}}

\begin{document}

\title[Black holes and class groups]{Black holes and class groups}

\author{Nathan Benjamin$^1$, Shamit Kachru$^1$, Ken Ono$^2$, and Larry Rolen$^3$}

\address{$^1$Stanford Institute for Theoretical Physics, Stanford University, Palo Alto, CA 94305}
\email{nathansb@stanford.edu, skachru@stanford.edu}

\address{$^2$Department of Mathematics and Computer Science, Emory University, Atlanta, GA 30322}
\email{ken.ono@emory.edu}

\address{$^3$Department of Mathematics, Vanderbilt University, 1326 Stevenson Center, Nashville, TN 37240}
\email{larry.rolen@vanderbilt.edu}

\begin{abstract}
The theory of quadratic forms and class numbers has connections to many classical problems in number theory.
Recently, class numbers have appeared in the study of black holes in string theory.
We describe this connection and raise questions in the hope of inspiring new collaborations between number theorists and physicists.
\end{abstract}

\dedicatory{In celebration of Don Zagier's 66th birthday}

\thanks{N.B. acknowledges support from the NSF and a Stanford Graduate Fellowship, S.K. is grateful for support from the NSF and a Simons Investigator Award, and K.O. thanks the Asa Griggs Candler Fund and the NSF for support. The authors thank Michael Mertens for helpful comments on an earlier version of this manuscript.}
\keywords{Black holes, class numbers, string theory}
\thanks{2010
 Mathematics Subject Classification: 11E41, 11R29,  83C57}

\maketitle

%%%%%%%%%%%%%%%%%%%%%%%%%%%%%%%%%%%%%%%%%%%%%%%%%%%%%%%%%%%%
\section{Introduction} \label{section1}
%%%%%%%%%%%%%%%%%%%%%%%%%%%%%%%%%%%%%%%%%%%%%%%%%%%%%%%%%%%%

In this note we will explain an interesting phenomenon which occurs at the interface of mathematics and physics.
 This was recently identified at a series of conferences on ``Number Theory, Geometry, Moonshine \& Strings'' at the Simons Foundation.
We aim to highlight the involvement in physics of objects which are well-known to number theorists, and which have arisen in new and somewhat surprising contexts in theoretical physics.
After literally centuries of exploration of the class numbers and class groups which are associated with binary quadratic forms of negative discriminant in the number theory literature, it has emerged recently that they appear to play an important (if somewhat mysterious) role in the physics of the simplest supersymmetric black holes in string theory.
The goal of this exposition -- aimed at both physicists and number theorists -- is to foster further research collaboration between mathematicians and physicists. For instance, we shall describe key facts about these objects from the number theory point of view, and ask whether there are corresponding physical interpretations. 

\medskip
As we mentioned, this new connection arises from considering {\it class numbers}, whose study goes back to Lagrange and Gauss. In particular, Gauss studied positive definite binary quadratic forms, namely polynomials of the form $[a,b,c](x,y):=ax^2+bxy+cy^2$ for fixed integers $a,b,c$. Two quadratic forms  $[a,b,c]$ and $[a',b',c']$ are said to be {\it equivalent} if there are integers $\alpha,\beta,\gamma,\delta$ with $\alpha\delta-\beta\gamma=1$ and $[a,b,c](x,y)=[a',b',c'](\alpha x+\beta y,\gamma x+\delta y)$. Under this relation, there are finitely many equivalence classes of forms of any fixed {\it discriminant} $D=b^2-4ac$, and this number is what is known as the class number of discriminant $D$. These equivalence classes also combine under an important operation known as the ``Gauss composition law,'' which turns them into a finite Abelian group $C(D)$. Such groups can also be understood  via ideal class groups of imaginary quadratic fields. See Section~\ref{section3} for more details on binary quadratic forms and class groups. Throughout this paper, by a {\it fundamental discriminant}, we mean a discriminant of an imaginary quadratic field.
Specifically, this means that $D\equiv0,1\pmod 4$, and that $D$ is square-free if $D\equiv1\pmod4$ and $D/4$ is a square free number congruent to $2$ or $3$ modulo $4$.

\medskip
From the theoretical physics side, the objects involved in this connection are \emph{black holes}.  These are solutions to Einstein's equations, first discovered by Schwarzschild in 1915.  Famously, they possess an interesting causal structure, including an ``event horizon'' which hides events in its interior from outside observers.  In the early 1970s, it was conjectured by Bekenstein that a black hole has a thermodynamic entropy given by the area of its event horizon (evaluated in units of the Planck area) \cite{Bekenstein}.  
This was confirmed, and elucidated considerably, in work by Bardeen, Carter and Hawking in classical general relativity \cite{Bardeen}, and in much further recent work studying the quantum statistical mechanics of black holes, where the entropy is given a statistical interpretation as a count of suitable microstates \cite{Strominger,Senlect}.
The connection to the work of Lagrange and Gauss arises because, as shall see in Section \ref{section2}, in
one of the simplest solutions of string theory (for a self-contained introduction to string theory, see, e.g.,
\cite{BBS}), it was found by Moore that the supersymmetric or BPS black hole solutions are in correspondence with the equivalence classes of binary quadratic forms 
introduced above.\footnote{By supersymmetric or BPS black holes, we mean those black hole solutions which preserve some of the supersymmetry of the underlying string theory. In general, states which preserve some supersymmetry are called BPS states.} The discriminant governs the black hole entropy, and the inequivalent classes at fixed discriminant correspond 
to distinct black hole solutions with the same entropy.
The particular objects Moore studied are the so-called supersymmetric attractor black holes of type IIB string theory compactified
on $K3 \times T^2$.  He classified these objects up to the duality symmetries of string theory, which provide non-perturbative
equivalences between naively distinct solutions, and are called ``U-dualities."

\medskip
We see, then, that the work of Moore, thanks to Gauss' group theoretic framework for binary quadratic forms,
may be arithmetically interpreted as the following theorem:

\begin{theorem}\label{Groups}
If $D<0$ is a fundamental discriminant, then the $U$-duality equivalence classes of attractor black holes of entropy $S = \pi \sqrt{-D}$ admit a structure as
 a finite Abelian group. Moreover, this group is isomorphic to the ideal class group of the field $\Q(\sqrt{D})$.
\end{theorem}

\begin{questions}
\begin{enumerate}
\item[]
\item[a)]
Is there a natural physical interpretation of the group law described in Theorem~\ref{Groups} in terms of attractor black holes? 
\item[b)] Is there a distinguished physical property of the identity class black hole, which corresponds to the class represented by the identity element $I_D$ given in \eqref{IDDefn}?
\item[c)] In there a physical relationship between inverse black holes (see \eqref{InversesQF})?
\item[d)]
What is the physical interpretation of the order of a black hole in the class group?
\end{enumerate}
\end{questions}

\medskip
These are just some of the questions that naturally come to mind in view of this correspondence between ideal classes and attractor black holes. In the last section of this paper, we offer a detailed discussion of possible future work. For completeness, before this discussion, we recall the background material which offers the connection between class numbers and black holes in Section~\ref{section2},
and in Section~\ref{section3} we recall classical work of Gauss on class groups of positive definite binary quadratic forms
and their relation to ideal class groups. 

%%%%%%%%%%%%%%%%%%%%%%%%%%%%%%%%%%%%%%%%%%%%%%%%%%%%%%%%%%%%
\section{Black holes and class numbers} \label{section2}
%%%%%%%%%%%%%%%%%%%%%%%%%%%%%%%%%%%%%%%%%%%%%%%%%%%%%%%%%%%%

The connection we explore here arose first in work of Moore, thoroughly explained in \cite{Moore1,Moore2}, with a telegraphic summary appearing in \cite{Moore3}.  These papers
studied prototypical examples of the {\it attractor mechanism} \cite{Ferrara} governing supersymmetric black holes in compactifications of string theory from ten to four dimensions.  
The string theory compactification involves (in the simplest case) a choice of Calabi-Yau threefold $X$, along with a point in its moduli space of Ricci-flat K\"ahler metrics.  In type IIB string theory, the resulting four-dimensional physical theory has $N=2$ supersymmetry (for a nice review of string compactification which focuses on 4d $N=2$ models and is accessible to mathematicians, see e.g. \cite{Aspinwall}).  The moduli of $X$ devolve into scalars in the vector multiplets (coming from variations of the complex structure of $X$) and the hypermultiplets (related to variations of the K\"ahler structure of $X$) of the $N=2$ supersymmetry.  The vector multiplets (together with a single graviphoton in the supergravity multiplet) also give rise to an Abelian gauge group of
rank $h^{2,1}(X) + 1$.

\medskip
In the setting of such string compactification, we wish to study black holes.  For ease, we study supersymmetry preserving black holes; this, and the ambient supersymmetry of the string compactification, can be viewed as furnishing a ``spherical cow" problem of the sort beloved by theoretical physicists.  Understanding the precise properties (such as solutions, entropies, and associated microstate counts) of such black holes is easier than for their non-supersymmetric counterparts, and so provides a natural starting point.

\medskip
The supersymmetric black holes in a 4d $N=2$ theory carry electric and magnetic charges under the $U(1)$ symmetries.
The charges are (loosely speaking) classified by a charge vector $Q \in H^{3}(X,{\mathbb Z})$.  For a given $Q$, the attractor mechanism
associates to the black hole with charge $Q$ a point $\tau_Q$ in the moduli space of complex structures on $X$ -- physically, this is the value the moduli approach at the horizon of the associated supersymmetric black hole \cite{Ferrara}. $\tau_Q$ is specified by
the fact that at the attractor point, the Hodge type of $Q$ becomes purely of type $(3,0) + (0,3)$:
\begin{equation}
Q \in H^{3,0}(X) \oplus H^{0,3}(X)~.
\label{blahblah}
\end{equation}
By counting, this involves eliminating $h^{2,1}(X)$ possible subspaces which could appear in the Hodge decomposition of $Q$, and so occurs at (complex) codimension $h^{2,1}(X)$ in the moduli space of complex structures -- i.e., at isolated points.

\subsection{Attractor black holes in $K3 \times T^2$ compactification}

Understanding the BPS black holes in a generic Calabi-Yau compactification is a rich problem that remains beyond reach.
However, for the case of $X = K3 \times T^2$, which enjoys reduced holonomy and consequently enhanced supersymmetry, 
one can characterize the BPS solutions completely.
The observation of Moore was that in this prototypical example, one can canonically associate the
attractor points to (equivalence classes of) binary quadratic forms of negative discriminant.  The discriminant $D < 0$ determines the associated attractor black hole entropy via
\begin{equation}
S = \pi \sqrt{-D}~.
\label{entropyis}
\end{equation}
Up to U-duality symmetries of the string compactification, there are a finite set of black holes at a given value of $-D$;
the class group acts to permute these.

\medskip
In greater detail, the moduli space of the string theory compactification on $K3 \times T^2$ is of the form
\begin{equation}\left(SO(22) \times SO(6)\right) \backslash SO(22,6;{\mathbb R}) / SO(22,6;{\mathbb Z}) \times  ({\bf H} / SL(2,{\mathbb Z}))~.\end{equation}
To unpack this slightly, we note that the moduli space of Ricci flat metrics on a K3 surface is given by
\begin{equation}\left(SO(19) \times SO(3) \right) \backslash SO(19,3;{\mathbb R}) / SO(19,3;{\mathbb Z}) \times {\mathbb R}^+~.\end{equation}
The additional dimensions present in the first factor above
(enhancing the (19,3) to (22,6)) correspond to scalars arising in dimensional reduction of the string theory antisymmetric tensor field, and from periods of the so-called ``Ramond-Ramond" (RR) higher-form gauge fields (which couple to D-branes) on the compact space.
One can think of the second factor ${\bf H} / SL(2,{\mathbb Z})$ as parametrizing the complex structure of the $T^2$.  The atypical structure of
this moduli space (as compared to generic Calabi-Yau threefolds, whose moduli locally factorize into products of vector multiplet and hypermultiplet moduli) is due to the extended 4d $N=4$ supersymmetry enjoyed by the $K3 \times T^2$ model.
Importantly, the bosonic fields in the 4d low-energy theory include 28 Abelian vector fields in addition to the (scalar) moduli.

\medskip
To find BPS states charged under these gauge groups, we need to wrap suitable D-branes on cycles of the compactification space.  To
understand this, note that the Abelian gauge fields arise from dimensional reduction of RR gauge fields under which D-branes
carry a charge \cite{Aspinwall, BBS}.  Therefore, wrapped D-branes are the particles which couple to the Abelian gauge fields in the 4d theory.
Of particular interest to us,
the IIB string theory has D3-branes that one can wrap on 3-cycles in the $X = K3 \times T^2$ to obtain supersymmetric black holes in the 4d low energy theory.  A general class in $H^3(X,{\mathbb Z})$ can be described as follows.  Choose a basis $\omega_i$, $i=1, \cdots, 22$ for $H^2(K3,{\mathbb Z})$.  Let $\alpha, \beta$ furnish a basis for $H^1(T^2,{\mathbb Z})$.  Then a basis for $H^3(X,{\mathbb Z})$ is
given by
\begin{equation}\omega_i \wedge \alpha, ~\omega_i \wedge \beta~.\end{equation}
We call the homology classes dual to the first set of forms the A-cycles, and those dual to the second set the B-cycles. The components of the charge vector along the A-cycles characterize the electric charge $E$ of the
black hole, while the components along the B-cycles give the magnetic charges $M$:
\begin{equation} E = \sum_i q_i \omega_i \wedge \alpha \equiv q \wedge \alpha, ~~M = \sum_i p_i \omega_i \wedge \beta \equiv p \wedge \beta ~. \end{equation}
One can find the attractors by expanding a general charge vector $Q \in H^3(X,{\mathbb Z})$ in this basis, and imposing the 
condition (\ref{blahblah}).  By explicit calculation together with a theorem of Shioda and Inose
\cite{ShiodaInose}, Moore finds that the attractor
points in the $K3 \times T^2$ moduli space combine:

\medskip
\noindent
$\bullet$
A choice of {\it singular} K3, i.e. a K3 surface of Picard rank 20.  The theorem of \cite{ShiodaInose} associates to such a K3 surface an elliptic curve.

\medskip
\noindent
$\bullet$ A value of the $\tau$ modulus of the compactification $T^2$ which determines an elliptic curve with complex multiplication; i.e. $\tau$ 
satisfies an equation of the form
\begin{equation}
a\tau^2 + b\tau + c = 0
\label{CM}
\end{equation}
with integer $a,b,c$.

\medskip
To a singular K3 is associated a natural two-dimensional lattice, the {\it transcendental lattice} $T$ (the complement in $H_2$ to the lattice spanned by the algebraic curve classes).  The attractor equations
tell us that we can take $p, q$ as generators of $T$.  Moore argues that the naturally associated quadratic form
\begin{equation}\begin{pmatrix} p^2 & -p \cdot q \\ -p \cdot q & q^2 
\end{pmatrix}\end{equation}
then has the same discriminant as (\ref{CM}), and the singular K3 surface is associated with the same $\tau$ value by the theorem
of Shioda and Inose.

\medskip
The discriminant
\begin{equation}
D = 4\left((p \cdot q)^2 - p^2 q^2 \right)< 0
\label{Dis}
\end{equation}
is invariant under the U-duality group, the $SO(22,6;{\mathbb Z}) \times SL(2,{\mathbb Z})$ 
arithmetic group which acts on the Teichmuller space.
The charges inequivalent under this group that give rise to the
same value of $-D$ are the distinct attractor black holes
with the same value of the black hole entropy, in the supergravity approximation.

\medskip
It follows from what we have said that there is a theorem:
\begin{theorem}[Moore]\label{UDualityThm}
The set of $U$-duality equivalence classes of attractor black holes with entropy $S = \pi\sqrt{-D}$ are in one-to-one correspondence with equivalence classes of primitive, positive definite binary quadratic forms of discriminant $D$. 
\end{theorem}

\medskip
An immediate corollary of this theorem is that there is a connection between counting functions for inequivalent attractor black holes at fixed entropy, and automorphic forms with coefficients governed by the class numbers of binary quadratic forms.  For an elementary discussion, see \cite{KachruTripathy}.

\begin{question} There are simple variants of the $K3 \times T^2$ compactification known as ``CHL strings," which involve other compactifications on complex threefolds of SU(2) holonomy.  What is the generalization of this story to CHL strings?
\end{question}

\begin{question} Does this phenomenon generalize in an interesting way to Calabi-Yau threefolds with generic holonomy? (Some conjectures regarding this question, and their preliminary exploration, may be found in \cite{Moore1,Moore2,Moore3}.) 
\end{question}

\medskip
So far, we only see the emergence of class numbers (as numbers of attractor points at a fixed entropy), and not the group structure
assigned to binary quadratic forms by Gauss.
We comment here on a conjectural physical interpretation of the class group.  While the U-duality group of the string compactification under discussion is $SO(22,6;{\mathbb Z}) \times SL(2,{\mathbb Z})$, the supergravity is ignorant of charge quantization.  Therefore supergravity naturally sees an $SO(22,6;{\mathbb R}) \times SL(2,{\mathbb R})$ symmetry acting on the moduli, and on solutions of the equations of motion.  As discussed in earlier papers by Sen \cite{Sengen} and particularly Cvetic and Youm \cite{CveticBPS}, one
can in fact generate all of the solutions with a given value of the black hole mass (controlled by $-D$), by starting with one
such solution and acting with suitable
$SO(22,6;{\mathbb R}) \times SL(2,{\mathbb R})$ ``solution-generating" transformations.  Of course most elements of the group give
rise to objects which do not exist in the full theory, as the charges are inappropriately quantized.  But for special choices, one
can obtain -- starting with a solution with correctly quantized charges -- other solutions with properly quantized choices and 
the same value of $-D$. 

\begin{question} In some (and perhaps all) cases, can we naturally view the
class group actions as arising from such $SO(22,6;{\mathbb R}) \times SL(2,{\mathbb R})$ symmetries of supergravity, generating the
orbit of $C(D)$ starting with a given solution?
\end{question}

\begin{question} A more general solution-generating group $O(24,8;\mathbb R)$, enlarging $SO(22,6;{\mathbb R}) \times SL(2,{\mathbb R})$, acts on static
solutions of the supergravity \cite{Sengen}.  Can this be used to find connections between BPS black holes at distinct values of 
$D$?  What is the mathematical interpretation?
\end{question}

\medskip
A further interesting correspondence is noted in \cite{Moore2}.
Suppose we label the attractor values of $\tau$ arising at a given $D$ as $\tau_i$, $i = 1, \cdots, h(D)$.  
Define $K_D = {\mathbb Q}(\sqrt{D})$, and its extension
by $j(\tau_i)$ to be $\hat K_D = K_D(j(\tau_i))$ (with $j$ the Klein $j$-function).  It turns out that this is independent of $i$ (by the theorem in \S3.5.4 of
\cite{Moore2}).
Then the attractor
$K3 \times T^2$ associated to $\tau_i$ is arithmetic, and is defined over (some finite
extension of) $\hat K_D$.  Moreover, the Galois group
$\Gal({{\hat K_D}/K_D})$ is isomorphic to the class group
$C(D)$, and acts to permute the attractors at a given value of 
$D$.  

\begin{question}
What is the connection between the supergravity duality group
and the Galois group $\Gal({{\hat K_D}/K_D})$?
\end{question}

\subsection{Examples} Here, we give some examples of possible embeddings of the class group in the solution-generating symmetry group. 
%We first describe an infinite class of embeddings of subgroups of class groups, and then specify to two simple explicit examples.
We first discuss two simple explicit examples mentioned in \cite{Moore1}, and then comment on some generalities below. 

\begin{example} $D=-20$:

\medskip
At $D=-20$, the class group is $\mathbb{Z}_2$ and there are two U-duality inequivalent black holes. One representative is for them to have charges
\begin{align}
p_1^2 = 1, ~~p_1\cdot q_1 &= 0, ~~q_1^2 = 5 \nonumber \\ p_2^2=2,~~p_2\cdot q_2&=1,~~q_2^2=3. 
\end{align}
Here we are using the conventions of \cite{CveticBPS, Sengen}.  In particular, of the 28 + 28 possible electric and magnetic charges, we work with 12 electric and 12 magnetic charges which rotate into one another under an $SO(6,6;{\mathbb R}) \subset SO(22,6;{\mathbb R})$ of the supergravity duality group.  These break into two sextuplets of electric and magnetic charges, with inner product governed by
the matrix (as in \cite{Sengen})

\begin{equation}
L=\begin{pmatrix} {\bf 0} & {\bf 1}\\ {\bf 1} & {\bf 0} \end{pmatrix}
\end{equation}
where ${\bf 0}$ and ${\bf 1}$ are $6\times 6$ matrices. An example choice of vectors $p_1, q_1, p_2, q_2$ above are
\begin{align}
p_1 &= \frac1{\sqrt{2}} \begin{pmatrix} 1 & 0 & 0 & 0 & 0 & 0 & 1 & 0 & 0 &0 &0 &0 \end{pmatrix} \nonumber \\  
q_1 &= \frac1{\sqrt{2}} \begin{pmatrix} 0 & 5 & 0 & 0 & 0 & 0 & 0 & 1 &0 &0 & 0 & 0\end{pmatrix} \nonumber \\
p_2 &= \frac1{\sqrt{2}} \begin{pmatrix} 2 & 0 & 0 & 0 & 0 & 0 & 1 & 0 & 0 &0 &0 &0\end{pmatrix} \nonumber \\
q_2 &= \frac1{\sqrt{2}} \begin{pmatrix} 0 &3  & 0 & 0 & 0 & 0 & 1 & 1 & 0 &0 &0 &0\end{pmatrix}.  
\end{align}
Acting with $O(6, 6, \mathbb{R}) \times SL(2,\mathbb{R})$ on $p_1$, $q_1$ to give $p_2$, $q_2$ gives us the constraints:
\begin{align}
\Omega(a p_1 + b q_1)^T &= p_2 \nonumber \\
\Omega(c p_1 + d q_1)^T &= q_2 \nonumber \\
a d - b c &= 1 \nonumber \\
\Omega^T L \Omega &= L. 
\end{align}
In order to match the T-duality invariants of $p^2$, $p\cdot q$, and $q^2$, we can choose the $SL(2,\mathbb R)$ element to be $\begin{pmatrix} -\frac{1}{\sqrt{23}} & \frac{3}{\sqrt{23}} \\ -\frac{8}{\sqrt{23}} & \frac{1}{\sqrt{23}} \end{pmatrix}$; this squares to negative the identity matrix. 

\end{example}
%We conjecture a solution to the above equations this where 
%\begin{align}
%\Omega^2 &\in O(6, 6, \mathbb{Z}) 
%\end{align}
%\end{example}

\begin{example} $D=-84$.

\medskip
At $D=-84$, the class group is $\mathbb{Z}_2\times \mathbb{Z}_2$, and there are four U-duality inequivalent black holes. One representative is for them to have charges
\begin{align}
p_1^2 = 1, ~~p_1\cdot q_1 &= 0, ~~q_1^2 = 21 \nonumber \\ p_2^2=3,~~p_2\cdot q_2&=0,~~q_2^2=7 \nonumber \\
p_3^2=2,~~p_3\cdot q_3&=1,~~q_3^2=11 \nonumber \\p_4^2=5,~~p_4\cdot q_4&=2,~~q_4^2=5.  
\end{align}
In the same basis as before, we can choose the following vectors:
\begin{align}
p_1 &= \frac1{\sqrt{2}} \begin{pmatrix} 1 & 0 & 0 & 0 & 0 & 0 & 1 & 0 & 0 &0 &0 &0 \end{pmatrix} \nonumber \\  
q_1 &= \frac1{\sqrt{2}} \begin{pmatrix} 0 & 21 & 0 & 0 & 0 & 0 & 0 & 1 &0 &0 & 0 & 0\end{pmatrix} \nonumber \\
p_2 &= \frac1{\sqrt{2}} \begin{pmatrix} 3 & 0 & 0 & 0 & 0 & 0 & 1 & 0 & 0 &0 &0 &0\end{pmatrix} \nonumber \\
q_2 &= \frac1{\sqrt{2}} \begin{pmatrix} 0 &7  & 0 & 0 & 0 & 0 & 0 & 1 & 0 &0 &0 &0\end{pmatrix} \nonumber \\
p_3 &= \frac1{\sqrt{2}} \begin{pmatrix} 2 & 0 & 0 & 0 & 0 & 0 & 1 & 0 & 0 &0 &0 &0 \end{pmatrix} \nonumber \\  
q_3 &= \frac1{\sqrt{2}} \begin{pmatrix} 2& 11 & 0 & 0 & 0 & 0 & 0 & 1 &0 &0 & 0 & 0\end{pmatrix} \nonumber \\
p_4 &= \frac1{\sqrt{2}} \begin{pmatrix} 5 & 0 & 0 & 0 & 0 & 0 & 1 & 0 & 0 &0 &0 &0\end{pmatrix} \nonumber \\
q_4 &= \frac1{\sqrt{2}} \begin{pmatrix} 4 & 5  & 0 & 0 & 0 & 0 & 0 & 1 & 0 &0 &0 &0\end{pmatrix}. 
\end{align}
It would be interesting if there are supergravity duality transformations acting on these vectors that realize $C(84)$.

\end{example}

Any black hole which has some charge representative with $p\cdot q=0$ must correspond to a class group element that squares to the identity. Suppose we have
\begin{align}
p^2 = a, ~~ p\cdot q = 0,~~ q^2 = -\frac{D}{4a}. 
\end{align}
We can choose the explicit vectors
\begin{align}
p &= \frac1{\sqrt{2}} \begin{pmatrix} 1 & 0 & 0 & 0 & 0 & 0 & a & 0 & 0 &0 &0 &0 \end{pmatrix} \nonumber \\  
q &= \frac1{\sqrt{2}} \begin{pmatrix} 0 & -\frac{D}{4a} & 0 & 0 & 0 & 0 & 0 & 1 &0 &0 & 0 & 0\end{pmatrix} 
\label{eq:asdf}
\end{align}
We can show that the following $SL(2,\mathbb R)$ element takes the ``canonical" identity electric and magnetic charges (i.e. with $p_{\text{Id}}^2=1$, $q_{\text{Id}}^2=-\frac{D}{4}$, $p_{\text{Id}}\cdot q_{\text{Id}} =0$) to have the same $T$-duality invariants as $p$ and $q$ in (\ref{eq:asdf}):
\begin{equation}
\begin{pmatrix} 0 & 2\sqrt{\frac{a}{-D}} \\   -\frac12\sqrt{\frac{-D}{a}} & 0 \end{pmatrix}.
\end{equation}
This element squares to minus the identity in $SL(2, \mathbb Z)$.
It would be interesting to show that there always is an $\Omega\in O(6,6;\mathbb{R})$ that acts on the identity charge vectors to form $p$, $q$, and whose square is in $O(6,6;\mathbb{Z})$.

%%%%%%%%%%%%%%%%%%%%%%%%%%%%%%%%%%%%%%%%%%%%%%%%%%%%%%%%%%%%
\section{Gauss's theory of class groups and Theorem~\ref{Groups}} \label{section3}
%%%%%%%%%%%%%%%%%%%%%%%%%%%%%%%%%%%%%%%%%%%%%%%%%%%%%%%%%%%%
In this section, we will describe the basic features of Gauss's composition law for positive definite binary quadratic forms. We will also briefly describe ideal class groups and how these relate to quadratic forms in the case of imaginary quadratic fields. More details about the topics discussed here can be found in \cite{Cox}, for instance, after which some of the following exposition is modeled. 
Throughout, $D$ denotes a negative, fundamental discriminant.  

\subsection{Gauss's composition law}

Gauss's composition law gives a procedure for combining two binary quadratic forms of discriminant $D$ and obtaining a new one, also of discriminant $D$. This law has its origins in the classical identity of Diophantus and Brahmagupta, which states that 
\[
(x^2+ny^2)(z^2+nw^2)=(xz+nyw)^2+n(xw-yz)^2. 
\]
In the language of quadratic forms given above, this can be interpreted as saying that if we multiply two values of the quadratic form $[1,0,n]$ together, we obtain another value of the same quadratic form. Gauss generalized algebraic identities such as these to study what he referred to as the {\it composition} of two forms (all forms of discriminant $D$): a composition of $[a,b,c]$ and $[a',b',c']$ is a form $[A,B,C]$ such that $[a,b,c](x,y)\cdot[a',b',c'](z,w)$ is always a value of $[A,B,C]$ at bilinear forms in the inputs $x,y,z,w$. To make the set of equivalence classes of binary quadratic forms into a group, we need an algebraic ``composition law'' which is well-defined as a function on classes, not just forms. Gauss describes a procedure for doing this, although his method is complicated. 

\medskip
Dirichlet gave a simpler method for obtaining such a composition law, which we describe here. 
\begin{definition}
Given binary quadratic forms $[a,b,c]$ and $[a',b',c']$ of discriminant $D<0$ with $\gcd(a,a',(b+b')/2)=1$, we define their {\bf composition} $[a,b,c]\circ[a',b',c']=[A,B,C]$ as the form with 
\[
A=aa',\quad\quad B=N, \quad\quad C=\frac{N^2-D}{4aa'}
,
\]
where $N$ is the unique integer modulo $2aa'$ such that 
\[
N\equiv b\pmod{2a},\quad\quad N\equiv b'\pmod{2a'},\quad\quad N^2\equiv D\pmod{4aa'}
.
\]
\end{definition}
This composition law preserves discriminant, primitivity (having all coefficients relatively prime), and, importantly, equivalence class. Thus  composition descends to a binary operation on the set $C(D)$ of classes of primitive positive definite binary quadratic forms of discriminant $D$.
\begin{theorem}
This set $C(D)$, under composition, is a group, called the {\bf class group}. \end{theorem}

\noindent
The order of $C(D)$ is called the {\bf class number}, and is denoted by $h(D)$. 
\bigskip

We now describe of the basic structure of these groups.
The identity element is represented by 
\begin{equation}
\label{IDDefn}
I_D:=
\begin{cases} 
[1,0,-D/4] & \text { if } D\equiv0\pmod 4,
\\
[1,1,\frac{1-D}4] & \text{ if } D\equiv1\pmod 4.
\end{cases}
\end{equation}
It is also easy to take inverses of elements under the composition law; explicitly, the inverse of the class represented 
 by a form $[a,b,c]$ can be represented by  
\begin{equation}
\label{InversesQF}
[a,-b,c].
\end{equation}
There is one other particularly well-understand facet of the structure of class groups. Specifically, 
Gauss's genus theory famously determines the $2$-torsion subgroup of $C(D)$, which we denote by $C_2(D)$.
\begin{theorem}[cf. Theorem 39, Corollary 1 of \cite{FroehlichTaylor}]
\label{ClassGrp2Torsion}
 For a fundamental discriminant $D$, $C_2(D)$ is an elementary $2$-group of order $2^{g-1}$, where $g$ is the number of distinct prime divisors of $D$.
 \end{theorem}
\begin{remark}
The class number, $h(D)$ is famously known to be finite, which isn't a priori clear from the definition above. To see why it is true, we need to introduce another definition. We call a quadratic form $[a,b,c]$ {\bf reduced} if either $-a<b\leq a< c$ or $0\leq b\leq a=c$. It turns out that every quadratic form is uniquely equivalent to a reduced form. Given this, one can directly check that there are only finitely many $a,b,c$ of fixed discriminant $D$ satisfying these inequalities, and this gives us a fast way to compute representatives for each of the classes.
\end{remark}

\subsection{Ideal class groups of imaginary quadratic fields}

It turns out that the composition operation on quadratic forms encodes deep information about arithmetic in quadratic fields. Specifically, for $d<0$, consider the imaginary quadratic field $\Q(\sqrt{d}):=\{a+b\sqrt d \ :\, a,b\in\Q\}$. We may assume that $d$ is square-free, and the {\it discriminant} of this field is 
\[
D:=\begin{cases}
d & \text{ if }d\equiv 1\pmod 4,
\\
4d & \text{ if } d\equiv 2,3\pmod 4.
\end{cases}
\]
 In this field, we consider the set of algebraic integers 
\[
\mathcal O_D:=\left\{a+b z_D \ :\, a,b\in\Z \right\},
\]
where 
\[
z_D:=\begin{cases}
\frac{1+\sqrt d}2& \text{ if }d\equiv 1\pmod 4,
\\
\sqrt d & \text{ if } d\equiv 2,3\pmod 4.
\end{cases}
\]
Associated to the number field $\Q(\sqrt D)$ is the {\it class group} of $\mathcal O_D$, which is the quotient of the group of so-called fractional ideals by the group of principal ideals. It turns at this ``ideal class group'' is isomorphic to the class group described in the previous subsection.
\begin{theorem}[cf. Theorem 7.7 of \cite{Cox}]\label{ClassGpIso}
Let $D<0$ be a fundamental discriminant. Then the ideal class group is isomorphic to $C(D)$. Explicitly, an isomorphism is given by mapping the quadratic form $[a,b,c]$ to the ideal $(a,(-b+\sqrt D)/2)$ in $\mathcal O_D$.
\end{theorem}

\subsection{Proof of Theorem~\ref{Groups}}
The proof follows by Theorem~\ref{UDualityThm} and Theorem~\ref{ClassGpIso}.

\subsection{Examples}

Here we consider the examples for the fundamental discriminants $D=-4,-84$.
We recall from Section \ref{section2} that the (equivalence classes of) binary quadratic forms at these values of the discriminant correspond to black 
holes with electric and magnetic charge vectors $q$, $p$ satisfying 
$$4\left(q^2p^2 - (q\cdot p)^2\right) = 4 ~(\text{resp.}~ 84)~,$$
counted up to U-duality action.  There are precisely two (four) such U-duality inequivalent black holes in the compactification on
$K3 \times T^2$.  Explicit representatives of the charge vectors for these black holes were presented at the end of Section \ref{section2}.

\begin{example}
Suppose that  $D=-4$, which corresponds to studying the integer ring $\mathbb Z[i]$ of Gaussian integers inside of $\mathbb Q(i)$. We search for reduced forms satisfying the inequalities above. First, we see if there are any solutions to the system
\[
4=4ac-b^2,\quad \quad a>0, \quad\quad -a<b\leq a< c.
\]
But these imply that $4>4a^2-a^2=3a^2$, and hence that $a=1$ is the only possible choice. Then the only possibilities for $b$ are $b=0$ or $b=1$. If $b=0$, then $4=4c$ implies $c=1$, and if $b=1$, we obtain $5=4c$, which isn't solvable in $\Z$. Similarly, if $0\leq b\leq a=c$, then $4\geq 4a^2-a^2=3a^2$ and so the only possibility is $a=c=1$, and hence $b=0$. All together, these imply that there is one class of quadratic forms, represented by $[1,0,1](x,y)=x^2+y^2$, and hence that the class number is $h(-4)=1$. The class group structure is trivial in this case.

\end{example}

\begin{example}
Here we consider the case when $D=-84$.
By the same reasoning as in the previous example,  the inequalities defining reduced forms imply that $3a^2\leq 84$, and so $1\leq a\leq 5$. 
A quick check then yields that the solutions to the set of inequalities $-a<b\leq a< c$ which correspond to our primitive forms are
$[1,0,21],[2,2,11],[3,0,7]$, while solving with the set of inequalities $0\leq b\leq a=c$ gives the unique solution $[5,4,5]$. Thus, 
$
C(-84)=\left\{[1,0,21],[2,2,11],[3,0,7],[5,4,5]\right\},
$
and $h(-84)=4$. 
It is not difficult to show that $C(-84)\cong\mathbb Z/2\Z\times\Z/2\Z$.

\end{example}

%%%%%%%%%%%%%%%%%%%%%%%%%%%%%%%%%%%%%%%%%%%%%%%%%%%%%%%%%%%%
\section{Further discussion} \label{section4}
%%%%%%%%%%%%%%%%%%%%%%%%%%%%%%%%%%%%%%%%%%%%%%%%%%%%%%%%%%%%

\subsection{Class Number Estimates}

The study of class numbers $h(D)$ and how they grow with discriminant is a classical problem in number theory. Gauss famously conjectured that for negative discriminants $D$,  we have $h(D)\rightarrow \infty$ as $D\rightarrow-\infty$. That is, Gauss surmised that for any natural number $n$, there are finitely many discriminants with $h(D)=n$. In particular, he conjectured complete lists of such $D$ for some special cases. For example, it is especially interesting to consider those discriminants with class number one. In this case, Baker showed in \cite{Baker} that there are only finitely many class number such discriminants, and Heegner and Stark (cf. \cite{Heegner} and \cite{Stark}) showed that the fundamental discriminants with class number $1$ are precisely $D=-3,-4,-7,-8,-11,-19,-43,-67,-163$, and the non-fundamental discriminants  are $D=-12, -16, -27, -28$. Gauss's conjecture is known in generality by Heilbronn \cite{Heilbronn}, and these explicit lists are known at least for class numbers $n=1,\ldots,100$ (thanks to work of Watkins \cite{Watkins}). 
 
\medskip
It is natural to ask for an asymptotic refinement of Heilbronn's theorem on the class number problem. Namely, as $D$ grows, how fast does $h(D)$ grow? The year following Heilbronn's work, Siegel gave an answer to this question, proving in \cite{Siegel} that
if $\varepsilon>0$, then there is a corresponding $c(\varepsilon)>0$ such that for sufficiently large $|D|$ we have that
\[
h(D)\geq c(\varepsilon) |D|^{\frac12-\varepsilon}.
\]
A short proof of this fact can also be found in \cite{Goldfeld}. The nature of Siegel's proof is quite interesting.
He obtained this conclusion explicitly assuming the truth of the Generalized Riemann Hypothesis (GRH), and he obtained the conclusion
inexplicitly (i.e. there is no formula for $c(\varepsilon)$) if the GRH turns out to be false. Either way, this offers an ineffective proof,
one which cannot even prove the infamous class number 1 problem discussed above. At present, the best known lower
bound is due to Goldfeld, Gross, and Zagier \cite{Goldfeld2}
$$
h(D) > \frac{1}{7000}\left(\log |D|\right) \prod_{\substack{p\mid D\\ p\neq |D|}} \left(
1-\frac{[2\sqrt{p}]}{p+1}\right).
$$
Improving on this lower bound is generally considered to be one of the most important problems in analytic number theory.

\begin{question}
How can these problems and theorems be interpreted in physics? 
For example, are there physical heuristics suggesting that the number of $U$-duality equivalence classes grow weakly as $D\rightarrow +\infty$?
\end{question}

\begin{figure}[h!]
  \caption{Plot of the class numbers associated to $\mathbb Q({\sqrt{-D}})$ as a function of $D$; the growth predicted by Siegel's theorem is evident.  Figure taken from J. Polak \cite{Polak}.}
  \label{fig:polak}
  \centering
    \includegraphics[width=0.8\textwidth]{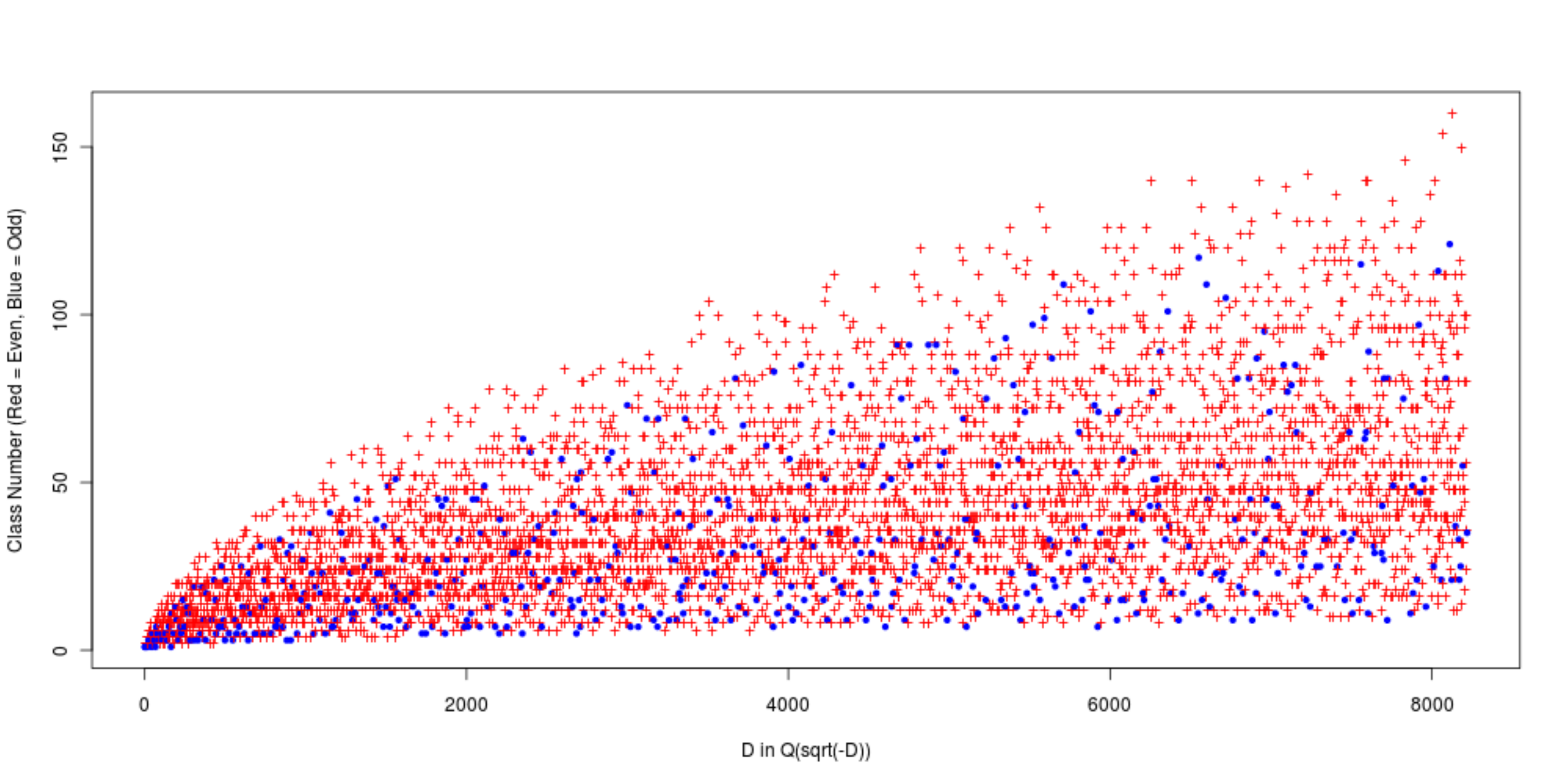}
\end{figure}

See Figure \ref{fig:polak} for a plot of class numbers.

\subsection{Class Number divisibility}

The subsection above was about the size of class groups, and hence numbers of equivalence classes of quadratic forms. In addition to studying their analytic growth rate, much work has been done on their arithmetic properties, such as the study of congruences they satisfy. For instance, in addition to studying the size of class groups, it is natural to ask about their structures as groups. For example, how large should we expect their $p$-torsion subgroups to be for primes $p$? For $p=2$, Gauss's genus theory gives a simple and elegant answer, as described in Theroem~\ref{ClassGrp2Torsion} above.

\medskip
Davenport and Heilbronn famously studied the case when $p=3$ (see \cite{DavenportHeilbronn}). In particular, they proved that the average size of the $3$-torsion subgroups of $h(D)$ tends to $2$, and that  the proportion of class numbers not divisible by $3$ is at least $1/2$. 
Little is known for $p\geq 5$.
However, known heuristics give conjectural answers towards understanding the structure of $p$-torsion subgroups of class groups. In particular, the expected probability that the $p$-torsion subgroup is trivial is predicted by the {\it Cohen-Lenstra {heuristics}} of \cite{CohenLenstra}, which in this special case claim that 
\[
\lim_{N\rightarrow\infty}\frac{\#\left\{-N<-D<0\ :\, p\nmid h(D)\right\}}{N}=\prod_{n\geq1}\left(1-\frac1{p^ n}\right).
\]

%More generally, the corresponding prediction for the probability that the $p$-rank is exactly $r$ is given by the expression
%\[
%p^{-r^2}\prod_{n\geq1}\left(1-p^{-n}\right)\prod_{i=1}^r\left(1-p^{-i}\right)^{-2}.
%\]
Other than the case when $p=3$ where substantial progress towards Cohen-Lenstra is known thanks to Davenport and Heilbronn, we are far from proving the predicted proportions. In fact, we are far from proving that a positive proportion of discriminants have class number not divisible by $p$. 
We do, however, have results such as Theorem 2 of \cite{KohnenOno},
which states that
\[
\lim_{N\rightarrow\infty}\frac{\#\left\{-N<-D<0\ :\, p\nmid h(D)\right\}}{N}\gg\frac{\sqrt{N}}{\log N}.
\]

Finally, we consider another structural question about class groups. For an integer $g\geq2$, let $N_g(x)$ {denote the number of square-free $D\leq x$} such that $C(-D)$ contains an element of order $g$. As discussed above, Gauss's genus theory predicts the answer for $g=2$. Specifically, the results above imply that a class group has an element of order $2$ if the discriminant has at least two prime factors, and so 
\[
N_2(x)\sim \frac{6x}{\pi^2}.
\]
For $g>2$, Cohen-Lenstra predicts that $N_g(x)\sim C_gx$, where
\[
C_g:=\frac{6}{\pi^2}\left(1-\prod_{i\geq1}\left(1-g^{-i}\right)\right).
\]
Towards this result, we have the following theorem of Soundararajan.
\begin{theorem}[Theorem 1 of \cite{Sound}]
We have that 
\[
N_g(x)\gg
\begin{cases}
x^{\frac12+\frac2g-\varepsilon}
& \text{ if } 
g\equiv0\pmod4
,
\\
x^{\frac12+\frac3{g-2}-\varepsilon}
& \text{ if } 
g\equiv2\pmod4
.
\end{cases}
.
\]
\end{theorem}

Since we do not have a specific speculation in mind, we highlight the following question.

\begin{question}
What does the structure of class groups mean in the context of attractor black holes? Can the physics interpretation of class groups put forth in this survey be used to study these old questions in number theory?
\end{question}

\subsection{Rademacher Sums}
Let $\Delta(q)$ be Ramanujan's Delta function, given explicitly by the formal power series expansion 
\[
\Delta(q)=q\prod_{n\geq1}\left(1-q^n\right)^{24}=:\sum_{n\geq1}\tau(n)q^n,
\]
where the Fourier coefficients are denoted by $\tau(n)$. As we shall see in the next section, the coefficients of $\Delta(q)$
are naturally expressed in terms of class numbers, and hence counts of $U$-duality equivalence classes of attractor
black holes. Therefore, we make a brief digression from class numbers and black holes to consider $\Delta(q)$, the prototypical
cuspidal Hecke eigenform in number theory.

\medskip
A distinguished physical role of not $\Delta(q)$, but rather $1/\Delta(q)$, is known. 
Counts of BPS states -- of interest for their connections to black hole physics or otherwise -- have been a major enterprise in theoretical
physics since the duality revolutions of the mid-1990s.  The simplest non-trivial state count -- of particle states preserving 1/2 of the supersymmetry in type IIA string compactification on K3, or the dual heterotic string compactification on $T^4$ -- gives as the
answer
\begin{equation}
\sum_{n} D_{n} q^{n-1} = {1\over \Delta(q)}~,
\end{equation}
where $D_n$ is interpreted as the number of BPS particle states of mass $n$.  This formula can be easily understood from the perspective
of the heterotic string: a 1/2 BPS state in heterotic string theory arises by keeping the right-moving string modes in their ground
state, and exciting the left-moving oscillators.  The left-moving degrees of freedom include (in light-cone gauge) 24 bosonic fields,
and the oscillator sums for these fields give the inverse eta-product characterizing $1/\Delta$.  It also has a dual description in
type II string theory: these BPS states arise from bound states of D0-branes to a D4-brane wrapping a K3 surface in IIA string theory, and one can understand $D_n$ as computing the Euler character of the Hilbert scheme of $n$ points on a K3 surface \cite{Vafa}:
\begin{equation}
\sum_n \chi({\rm Hilb}^n(K3)) q^{n-1} = {1\over \Delta(q)}~.
\end{equation}

\medskip
One useful way to think of these formulae is in terms of the Rademacher sums of number theory, {which are often called Poincar\'e series in the modular forms literature.}
For physicists, one can think of a Rademacher series in the following heuristic way.  Suppose one has a symmetry group
$G$, and one wishes to construct a function $f$ invariant under $G$ action.  One simple method is to construct $f$ by starting
with a single term -- say a constant -- and then summing over $G$ - images of this term.  As long as the sum converges, the
resulting function will be $G$-invariant by construction.

\medskip
Rademacher's method similarly constructs a modular form by starting with the most polar term in its $q$-expansion, and then successively summing over
$\Gamma = SL(2,{\mathbb Z})$ images of this term.  More properly, because $SL(2,{\mathbb Z})$ is infinite, one should sum over cosets
$\Gamma_{\infty} \backslash \Gamma$ (with $\Gamma_\infty$ defined as the upper triangular matrices in $SL(2,{\mathbb Z})$) in order to obtain a convergent expression. 
{Some care has to be taken with convergence of the sum when studying modular forms of negative weight.}     This formal construction obtains physical interpretation in the context of string theory, and 
especially in dualities between AdS$_3$ gravity and 2d conformal field theories.  There, in the ``Farey Tail" of Dijkgraaf et al. 
\cite{FareyTail}, the Rademacher sum representations of suitable partition functions are interpreted as sums over saddle points in 
AdS$_3$ quantum gravity, with different powers of $q$ weighing different saddles by an appropriate Euclidean action for the 
corresponding gravity solution.

\medskip
We can apply Rademacher's expansion to $1/\Delta$ (a case originally studied by Rademacher himself). 
As usual, we let $K(m,n;c)$ be the {\it Kloosterman sum}
\[
K(m,n;c):=\sum_{d\in(\Z/c\Z)^{\times}}e^{2\pi i \frac{m\bar d+nd}c},
\]
where $\bar d$ denotes  the multiplicative inverse of $d$ modulo $c$. 
Then it is a classical fact (for example, see \S 6.3 of \cite{Colloq}) that the Fourier coefficients of
\[
\frac{1}{\Delta(q)}=:q^{-1}+\sum_{n\geq 0}a(n)q^n,
\]
for $n\geq1$ are given by the Rademacher sums
\begin{equation}\label{RademacherSum1}
a(n)=\frac{2\pi}{n^{\frac{13}{2}}}\cdot \sum_{c>0}\frac{K(-1,n;c)}{c}\cdot I_{13}\left(\frac{4\pi\sqrt n}c\right),
\end{equation}
where $I_{13}$ is the index 13  $I$-Bessel function. 
This expansion actually should have a nice path integral sum-over-saddles interpretation in string theory as giving the generating function of entropies of 
``small black holes" in compactification on $AdS_2 \times S^2 \times K3 \times T^2$, though details remain hazy; see for instance \cite{Gomes}.

\medskip
In number theory, $\Delta(q)$ is more prominent than $1/\Delta(q)$. It is the prototype of a Hecke eigenform, and as a result
offered early examples of deep structures and results such as  Galois representations, the Weil conjectures, and Hecke theory. 
It is very interesting to note that both $1/\Delta(q)$ and $\Delta(q)$ are {unified} in that their Fourier expansions are both ``Rademacher sums''.
In the case of $\Delta(q)$,  for any $n\geq2$ we have that
\[
\tau(n)
=\frac{2\pi n^{\frac{11}2}}{\beta_{\Delta}}\cdot \sum_{c>0}\frac{K(1,n;c)}{c}\cdot J_{11}\left(\frac{4\pi\sqrt n}c\right)
\]
where $\beta_{\Delta}=2.840\dots$, and $J_{11}$ is the usual
order $11$ $J$-Bessel function.

\medskip
Since Rademacher sums appear naturally in physics, it is natural to ask the following question.
\begin{question}
Is there a physical interpretation of $\Delta(q)$ given its formulation in terms of ``holomorphic'' Rademacher sums?
\end{question}
It is immediately evident to every physicist that $\Delta$ is closely related to a torus partition function in a 2d CFT of
24 free fermions with a particular spin structure, but in fact that function vanishes due to fermion zero modes, and giving
any precise physical insight into facts about the Fourier coefficients $\tau(n)$ of $\Delta$ remains an intriguing open problem.

\medskip
We close this section by noting that the formulae for counting of 1/2 BPS states explored here have a rich generalization to
the counting of 1/4 BPS states in string compactification on $K3 \times T^2$, first explored in \cite{DVV}.  Connections of the 1/4 BPS state counts
with the theory of Siegel forms, Jacobi forms, and mock modular forms, and associated questions in number theory, were discussed in
the work of Dabholkar, Murthy and Zagier \cite{DMZ}.

\subsection{Eichler-Selberg Trace formula}

The function $\Delta(q)$ can also be interpreted as the infinite collection of Hecke eigenvalues, as the Fourier coefficients $\tau(n)$ are precisely the eigenvalues of $\Delta$ under the action of the $n$-the Hecke operator $T_n$ acting on the one-dimensional space of weight $12$ cusp forms for $\mathrm{SL}_2(\Z)$. 
{The well-known {\it Eichler-Selberg trace formula} gives an explicit formula for the trace of $T_n$ acting on the space of weight $k$ cusp forms on  $\mathrm{SL}_2(\Z)$ in terms of class numbers (see, e.g., \cite{Zagier}). }
To give this formula, we need to define a few notations. Firstly, we need the {\it Hurwitz class numbers} $H(D)$, a modification of the $h(D)$, which count the number of equivalence classes of positive definite binary quadratic forms (not necessarily primitive) of discriminant $D$, and in this count weights forms which are equivalent to a form of the shape $[a,0,a]$ by a $1/2$ instead of a $1$ and weights forms equivalent to a form of the shape $[a,a,a]$ by a $1/3$ (the ``$2$'' and ``$3$'' here count the order of the stabilizer of the form in $\mathrm{SL}_2(\Z)$, which is $1$ in all other cases).  Letting $\omega_D$ denote this stabilizer count, which is generically $1$ but is $2$ or $3$ in the distinguished cases just listed which correspond to $D=-4,-3$, respectively, we also have the explicit formula
\[
H(D)=\sum_{d^2|D}h\left(\frac{D}{d^2}\right)\omega_D^{-1}.
\]
Furthermore, we let $P_k(t,N)$ denote the coefficient of $x^{k-2}$ in the power series expansion of 
\[
\frac1{1-tx+Nx^2}.
\]
The Eichler-Selberg trace formula then says that the trace of the Hecke operator $T(n)$ on the space of weight $k$ cusp forms for any even weight $k\geq4$ is:
\[
-\frac12\sum_{t\in\Z}P_k(t,n)H(t^2-4n)-\frac12\sum_{dd'=n}\min(d,d')^{k-1}
,
\]
where we let $H(n)=0$ for $n>0$, and we formally define $H(0)=-\frac1{12}$.
When $k=12$, the space of cusp forms is one-dimensional, spanned by $\Delta$, and so the trace of $T(n)$ is simply the eigenvalue $\tau(n)$. 
For example, for a prime $p$, these formulas then say that 
\[
\tau(p)=-\frac12\sum_{t=0}^{\left\lfloor\sqrt{p}\right\rfloor}\left(t^{10}-9pt^8+28p^2t^6-35p^3t^4+15p^4t^2-p^5\right)H(t^2-4p)
-\frac12\sum_{dd'=p}\min(d,d')^{11}
.
\]

\begin{question}
Since the Eichler-Selberg traces formulas for $\Delta(q)$ are phrased in terms of class numbers, it is natural
to ask: do the coefficients of $\Delta(q)$ have a natural physical interpretation in terms of black holes?
\end{question}

\subsection{Counting elliptic curves}
A final interpretation of class numbers we discuss here realizes them as counting isomorphism classes of elliptic curves with fixed structure (based on seminal work of Deuring).
For example, if we specify the total number of points over a finite field $\mathbb F_q$ for a prime power $q=p^r$, then we can give explicit formulas for these isomorphism counts. 
For instance, we have the following sample result of Schoof, where in the following, $I(t)$ denotes the isogeny class of elliptic curves with $q+1-t$ points over $\F_q$ and $N(t)$ is the number of $\mathbb F_q$-isomorphism classes of curves in $I(t)$.
\begin{theorem}[Special case of Theorem 4.6 in \cite{Schoof}]\label{Schoof1}
If $t^2<4p$ and $p\nmid t$ or $t=0$ and $q$ is not a square, then 
\[
N(t)=H(t^2-4p).
\]
\end{theorem}
\begin{remark}
Schoof also gives formulas for $N(t)$ in other cases as well, but the $6$ other cases all have simple formulas which are boundary cases which don't involve class numbers, so we omit them.
\end{remark}
Note that in particular, Theorem~\ref{Schoof1} applied to the {\it supersingular case} over a prime-order field $\mathbb F_p$ with $t=0$ implies that $N(t)=H(-4p)$. 
Schoof also proves results {on the number} of isomorphism classes of elliptic curves over $\F_q$ with distinguished torsion subgroups, such as the following.
\begin{theorem}[Special case of Theorem 4.9 in \cite{Schoof}]
If $n\geq1$ is odd, $t^2\leq 4q$, $p\nmid, t$, $q\equiv 1\pmod n$, and $t\equiv q+1\pmod{n^2}$, then
\begin{align*}
&
\#\left\{
\F_q\text{-isomorphism classes of elliptic curves } E \text{ in } I(t)  :\  E(\F_q)[n]\cong\Z/n/Z\oplus/Z/n\Z
\right\}
\\
=
&
H\left(\frac{t^2-4q}{n^2}\right)
.
\end{align*}
\end{theorem}
 The proofs of these theorems identify ideal classes with elliptic curves through the theory of complex multiplication.
 Supersingular elliptic curves (i.e. those with $t=0$) stand out prominently in number theory. For example, the locus
 of such curves is a fundamental device in the theory of $p$-adic modular forms. In view of the special role of these curves, it is natural to ask the following question.

\begin{question}
{Considering $t=0$, i.e., in the case of supersingular elliptic curves in characteristic $p$, are there corresponding special properties of the attractor black holes with entropy $S = \pi \sqrt{4p}$?}
\end{question}

\subsection{Extremal conformal field theories}

In this subsection, we will discuss one further place that class numbers have made a brief appearance in physics. The AdS/CFT correspondence \cite{Juan} relates a $d$-dimensional conformal field theory (CFT) with a gravity theory in $d+1$-dimensional Anti de Sitter space. In \cite{Witten}, a natural question was asked: what two-dimensional CFTs would be dual to a putative theory of pure three-dimensional gravity, with only graviton excitations in the spectrum, along with black holes (known as ``BTZ black holes" in 3d gravity)?

\medskip
In \cite{Witten}, the assumption made was that the two-dimensional conformal field theory would have a spectrum ``as close as possible" to the vacuum character of the CFT. A large simplifying assumption was also made -- that the two-dimensional conformal field theory would holomorphically factorize, so that the mathematical machinery of modular forms could be used. An \emph{extremal CFT} was defined in \cite{Hohn1, Hohn2} to be a (hypothetical) holomorphic CFT at central charge $c=24k$ that had partition function
\begin{equation}
Z_k(q) = q^{-k} \prod_{n=2}^{\infty} \frac{1}{(1-q^n)} + \mathcal{O}(q^1)
\end{equation}
such that, at fixed integer $k$, $Z_k(q)$ was the unique function satisfying both the above constraint, and was modular invariant under $SL(2,\mathbb{Z})$ transformations. In \cite{Witten} an elementary argument was given showing that the function $Z_k(q)$ always exists, but it is unknown which extremal CFTs exist for $k>1$.\footnote{At $k=1$, the theory of Frenkel-Lepowsky-Meurman \cite{FLM} which plays a starring role in Monstrous moonshine is extremal.} Nevertheless, it turns out that these modular functions have a description that involves both Rademacher sums and the CM points counted by these class numbers. In the spirit of this paper, we record these identities here.

\medskip
For any positive integer $d$, we define the Rademacher sum 
$
R_d(\tau) = q^{-d} + \sum_{n\geq1} r_{d,n}q^n,
$
where
\[
r_{d,n}
= 2 \pi \sqrt{\frac{d}{n}}
 \times \sum_{c > 0} \frac{K(-d,n;c)}{c}
 \cdot
I_{1}\left(\frac{4\pi \sqrt{dn}}{c}\right),
\]
and $I_1$ is a modified $I$-Bessel function.
Furthermore, we recall a finite, algebraic formula for the partition numbers $p(n)$ obtained in \cite{BruinierOno2}.
This is expressed in terms of the quasimodular Eisenstein series $E_2$ and
Dedekind's eta function.
We then define the weight $-2$,
level $6$ modular function
\[
G(\tau):=\frac{1}{2}\frac{E_2(\tau)-2E_2(2\tau)-3E_2(3\tau)+6E_2(6\tau)}{\eta(\tau)^2\eta(2\tau)^2\eta(3\tau)^2\eta(6\tau)^2}
\]
and its non-holomorphic derivative 
\[
\mathcal{P}(\tau):=\frac{i}{2\pi}\frac{\partial G}{\partial\tau}-\frac{G(\tau)}{2\pi\operatorname{Im}(\tau)}.
\]

\medskip
{For positive integers $D$, let 
 \[
\mathcal{Q}_{D}:=
 \Gamma_0(6)\backslash\left\{Q=[a,b,c]: a,b,c\in\Z, b^2-4ac= -24D+1, 6|a, a>0,\ b\equiv 1 \pmod{12} \right\},
\]
where $\Gamma_0(6)$ acts on the set of quadratic forms in the usual way.
For any $Q\in \mathcal Q_D$, we also take $\tau_Q$ to be the CM point in the upper half plane satisfying the $a\tau_Q^2+b\tau_Q+c=0$.}
It is known that $(24n-1)P(\tau_Q)$ is an algebraic integer for any $Q\in\mathcal{Q}_{D}$.
Thus,  it is reasonable to consider  the
 traces of singular moduli
\[
\operatorname{Tr}(\mathcal P;n):=\sum_{Q\in\mathcal{Q}_{n}}\mathcal{P}\left(\tau_Q\right)
.
\]
In \cite{BruinierOno2}, it was shown that 
$
\operatorname{Tr}(\mathcal P;n)=(24n-1)p(n).
$

\medskip
By combining these facts, we obtain the following closed form expression of these proposed extremal partition functions in terms of sums over generalized class groups and Rademacher sums (see  \cite{WittenNote}):
\[
Z_k(q)=\frac{1}{24k-1}\operatorname{Tr}(\mathcal P;k)+\left(R_k(\tau)-R_{k-1}(\tau)\right)+\sum_{n=1}^{k-1}\frac{1}{24n-1}\operatorname{Tr}(\mathcal P;n)\left(R_{k-n}(\tau)-R_{k-n-1}(\tau)\right).
\]

\medskip
In \cite{GGKMO}, an extension of extremal CFTs to $N=(2,2)$ supersymmetry was defined. There, instead of considering purely holomorphic partition functions (as would be associated to chiral CFTs), we can instead consider the elliptic genus (a holomorphic index), and demand that the elliptic genus matches the $N=2$ superconformal vacuum descendants. The elliptic genus of an $N=(2,2)$ superconformal field theory with integral $U(1)$ charges can be shown to transform as a weak Jacobi form of weight $0$ and index $m=c/6$ (where $c$ is the central charge of the CFT) \cite{Kawai}.

\medskip
The analogous definition of an extremal $N=(2,2)$ CFT would be that the \emph{polar} pieces of the elliptic genus vanish. In particular, define the Fourier coefficients of the elliptic genus as follows:
\begin{equation}
Z_{\text{EG}}(q,y) = \sum_{n, \ell} c(n, \ell) q^n y^{\ell}
\end{equation}
and define the \emph{polarity} to be 
\begin{equation}
p(n, \ell) = 4mn-\ell^2
\end{equation}
where $m$ is the index of the weak Jacobi form (and related to the central charge of the CFT by a factor of $6$). Terms in the expansion of $Z_{\text{EG}}(q,y)$ with $p(n, \ell) < 0$ are called polar terms, and are interpreted as being dual to particle states, in contrast to nonpolar terms (which have $p(n,\ell)>0$) being interpreted as black hole states. An extremal $N=(2,2)$ theory would be one that matched the polar terms to the superVirasoro vacuum descendants. In \cite{GGKMO}, it is shown that except for a finite list of $m$ (namely, $m=1, 2, 3, 4, 5, 7, 8, 11$, and $13$), there does not exist a function that both transforms with the right modularity properties, and reduces to the superVirasoro descendants for all polar terms (see also \cite{Manschot}).
\footnote{For $m=2,4$, explicit CFTs which realize this exist, and are constructed in \cite{Cheng,Benjamin}.}
 This is to be contrasted with the nonsupersymmetric case in \cite{Witten}.  There, although it is unknown whether or not the CFTs exist, for all values of the central charge, there exists a candidate partition function with the right modular properties that reduces to the Virasoro vacuum descendants for the polar pieces.  

\medskip
A brief sketch of the argument supplied in \cite{GGKMO} goes as follows. It can be shown \cite{Eichler} that at index $m$, the dimension of the space of weak Jacobi forms, $J(m)$, and the number of independent polar terms, $P(m)$, grows differently with $m$. In particular, we have
\begin{equation}
J(m) = \left \lfloor{\frac{m^2}{12}+\frac{m}2+1}\right \rfloor 
\end{equation}
whereas
\begin{equation}
P(m) = \frac{m^2}{12} + \frac{5m}8 + \frac14 \sum_{d|4m} h(d) - \frac12 \left \lfloor \frac b2 \right \rfloor - \frac12 ((\frac m4)) + \frac1{24}
\end{equation}
where $h(d)$ is the class number associated with discriminant $-d$ (with the exception of the cases $d=3,4$, where $h(3) = 1/3$, and $h(4) = 1/2$ in this formula); $b$ is the largest integer such that $b^2 | m$; and 
\begin{equation}
((x)) = x - \frac12 \left( \left \lceil x \right \rceil + \left \lfloor x \right\rfloor\right).
\end{equation}
Since at large $m$, we have $P(m)$ growing larger than $J(m)$, in general we are unable to ``tune" all of the polar terms to match the superVirasoro vacuum character.

\begin{figure}[h!]
  \caption{A plot of $\frac{P(m) - \frac{m^2}{12} - \frac{5m}8}{\sqrt{m}}$ for the first $100000$ values of $m$.}
  \vspace*{0.2in}
  \label{fig:scatter}
  \centering
    \includegraphics[width=0.7\textwidth]{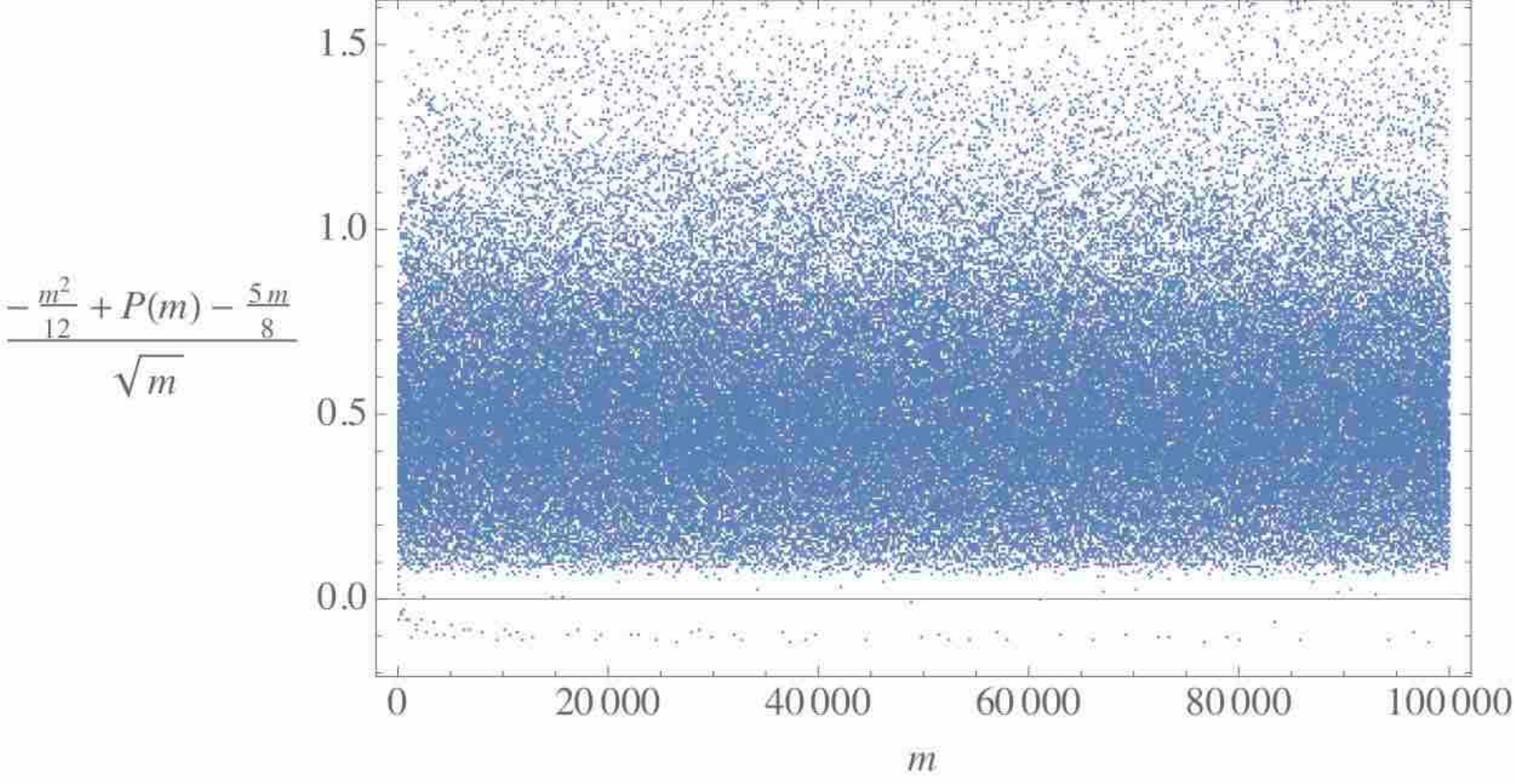}
\end{figure}

\begin{figure}[h!]
  \caption{Cumulative probability distribution for $\frac{P(m) - \frac{m^2}{12} - \frac{5m}8}{\sqrt{m}}$ calculated from the first $100000$ values of $m$.}
  
  \smallskip
  \label{fig:cumulative}
  \centering
    \includegraphics[width=0.7\textwidth]{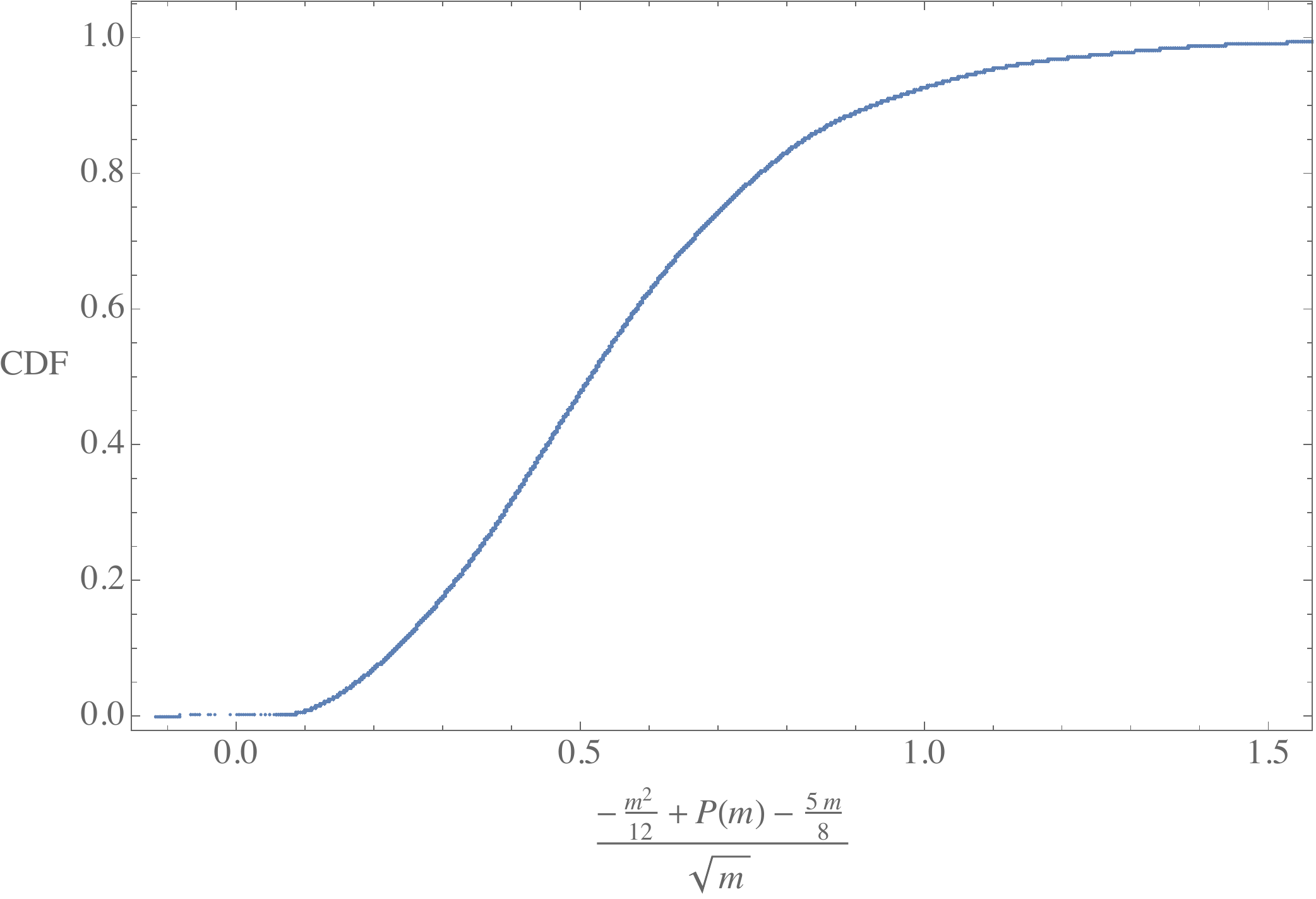}
\end{figure}

\begin{figure}[h!]
  \caption{Probability distribution for $\frac{P(m) - \frac{m^2}{12} - \frac{5m}8}{\sqrt{m}}$ calculated from the first $100000$ values of $m$.}
  \label{fig:prob}
  \centering
    \includegraphics[width=0.7\textwidth]{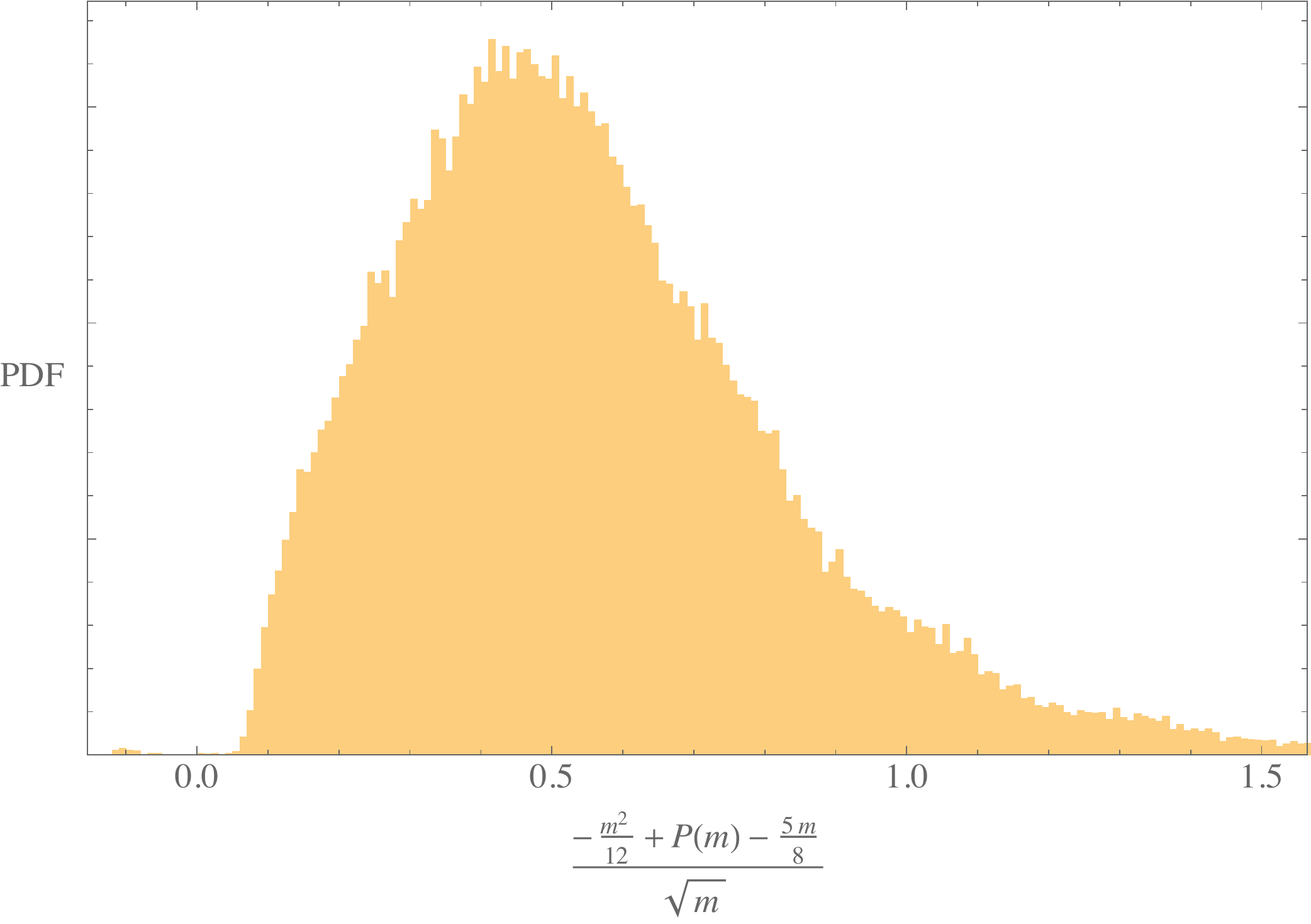}
\end{figure}

\medskip
In \cite{SmallBH}, a physical explanation was proposed for the $\mathcal O(m)$ term in the difference between $P(m)$ and $J(m)$ in terms of small black holes. However, the remaining terms were unexplained. Roughly speaking, the remaining terms grow as a random number multiplying $\sqrt{m}$ (see Fig. \ref{fig:scatter}, \ref{fig:cumulative}, \ref{fig:prob}); it would be interesting if we could get an understanding of this term in terms of small black holes as well.


%{\red Anyone know how to force these figures to appear in Section 4.6?} {\blue I added a new page command -- seems to work for now, let's make sure it's the right place before we resubmit}
 
\newpage

\subsection{Non-supersymmetric attractors}

The attractor mechanism has a generalization to non-BPS black holes as long as they are extremal, saturating a condition relating
their mass to their electric and magnetic charges \cite{Goldstein}.  Simple examples in conventional general relativity include the
Reissner-Nordstrom black holes.  This more general non-supersymmetric attractor mechanism will apply to extremal black holes in
theories with moduli spaces of vacua, for instance arising in supersymmetric theories as non-supersymmetric solutions.  

\medskip
The equations for both supersymmetric and non-supersymmetric attractors in IIB compactification on a Calabi-Yau threefold $X$ with
holomorphic three-form $\Omega$ can be described as follows.  Choose a charge $Q \in H^3(X,{\mathbb Z})$, and define
\begin{equation}Z_Q = \int_X Q \wedge \Omega~.\end{equation}
Let $K$ denote the K\"ahler potential on the complex structure moduli space of $X$.  Define the effective potential
\begin{equation}V_Q = e^K \left( \sum_{\alpha} |D_\alpha Z_Q|^2 + |Z_Q|^2 \right)~.\end{equation}
Here, the $\alpha$ index runs over a label on the complex structure moduli of $X$, and $D_\alpha Z_Q = \partial_{\alpha}Z_Q + K_{\alpha} Z_Q$.

\medskip
The BPS attractors on $X$ are precisely those which extremize both of the terms in $V_Q$ independently.  The non-supersymmetric
attractors on $X$ are instead those which minimize $V_Q$, but do not minimize each term separately \cite{Tripathy}.  Their physics, geometry, and
number theory has been much less studied than that of their BPS counterparts (where in turn, only the surface has been scratched).

\begin{question} Is there an interesting geometric or number theoretic significance to the non-supersymmetric attractor points in
the $K3 \times T^2$ moduli space, generalizing the connection of BPS black holes to class groups?
\end{question}


\begin{thebibliography}{99}


\bibitem{Aspinwall}   P.~S.~Aspinwall,
  \emph{Compactification, geometry and duality: N=2},
  arXiv:hep-th/0001001.
  %%CITATION = HEP-TH/0001001;%%

\bibitem{Baker} A. Baker, {\it Linear forms in the logarithms of algebraic numbers}, Mathematika {\bf 13} (1966), 204--216.

\bibitem{Bardeen} J. Bardeen, B. Carter and S. Hawking, \emph{The four laws of black hole mechanics}, Comm. Math. Phys. {\bf 31} (1973), 161--170.



\bibitem{BBS}  K.~Becker, M.~Becker and J.~H.~Schwarz,
  \emph{String theory and M-theory: A modern introduction}, Cambridge University Press (2007).

\bibitem{Bekenstein} J. Bekenstein, \emph{Black holes and entropy},
Phys. Rev. D{\bf 7} (1973), 2333--2346.

\bibitem{Benjamin} 
  N.~Benjamin, E.~Dyer, A.~L.~Fitzpatrick and S.~Kachru,
  \emph{An extremal ${\mathcal{N}}=2$ superconformal field theory},
  J.\ Phys.\ A {\bf 48}, no. 49 (2015), 495401 [arXiv:1507.00004 [hep-th]].


\bibitem{SmallBH} 
  N.~Benjamin, E.~Dyer, A.~L.~Fitzpatrick, A.~Maloney and E.~Perlmutter,
  \emph{Small Black Holes and Near-Extremal CFTs},
  JHEP {\bf 1608}, 023 (2016)
  [arXiv:1603.08524 [hep-th]].


  
\bibitem{Colloq}  K. Bringmann, A. Folsom, K. Ono, and L.Rolen,
\emph{Harmonic Maass forms and mock modular forms: Theory and Applications},
Amer. Math. Soc., 2017.

\bibitem{BruinierOno2} 
J. Bruinier and K. Ono,
\emph{Algebraic formulas for the coefficients of half-integral weight
  harmonic weak {M}aass forms,} Adv. Math. \textbf{246} (2013), 198--219.


\bibitem{Cheng} M.~C.~N.~Cheng, X.~Dong, J.~F.~R.~Duncan, S.~Harrison, S.~Kachru and T.~Wrase,
  \emph{Mock Modular Mathieu Moonshine Modules},
  Res.\ Math.\ Sci.\  {\bf 2} (2015), 13 
  [arXiv:1406.5502 [hep-th]].

\bibitem{CohenLenstra} H. Cohen and H.W. Lenstra, \emph{Heuristics on class groups of number fields}, Springer Lect. Notes in Math. {\bf 1068} (1984), 33--62.

\bibitem{Cox} D. Cox, \emph{ Primes of the Form $x^2 + ny^2$. Fermat, Class Field Theory and Complex Multiplication}, Wiley, 1989.

\bibitem{CveticBPS} M. Cvetic and D. Youm, \emph{Dyonic BPS saturated black holes of heterotic string on a six torus},
Phys. Rev. D{\bf 53} (1996), 584--588 [hep-th/9507090].

\bibitem{DMZ} A. Dabholkar, S. Murthy and D. Zagier, 
\emph{Quantum Black Holes, Wall Crossing, and Mock Modular
Forms}, arXiv:hep-th/1208.4074.

\bibitem{DavenportHeilbronn} H. Davenport and H. Heilbronn, \emph{On the density of discriminants of cubic fields II}, Proc. Roy. Soc. London Ser. A {\bf 322} no. 1551 (1971) 405--420.


\bibitem{FareyTail} R. Dijkgraaf, J. Maldacena, G. Moore and
E. Verlinde, \emph{A black hole Farey Tail}, arXiv:hep-th/005003.

\bibitem{DVV} R. Dijkgraaf, E. Verlinde and H. Verlinde,
\emph{Counting dyons in $N=4$ string theory}, Nucl. Phys. 
B {\bf 484}, 543 (1997).


\bibitem{Eichler} M. Eichler and D. Zagier, \emph{The theory of Jacobi forms}, Birkhauser (1985).

\bibitem{Ferrara} 
  S.~Ferrara, R.~Kallosh and A.~Strominger,
  \emph{N=2 extremal black holes},
  Phys.\ Rev.\ D {\bf 52} (1995), R5412 [hep-th/9508072].
 

\bibitem{FLM} I. Frenkel, J. Lepowsky and A. Meurman, \emph{Vertex Operator Algebras and the Monster},
Academic Press (1988).

\bibitem{FroehlichTaylor} A. Fr\"ohlich and M. Taylor, \emph{Algebraic number theory}, volume 27 of Cambridge Studies in Advanced Mathematics.
Cambridge University Press, Cambridge, 1993.

\bibitem{GGKMO} 
  M.~R.~Gaberdiel, S.~Gukov, C.~A.~Keller, G.~W.~Moore and H.~Ooguri,
  \emph{Extremal N=(2,2) 2D Conformal Field Theories and Constraints of Modularity},
  Commun.\ Num.\ Theor.\ Phys.\  {\bf 2}, 743 (2008)
  [arXiv:0805.4216 [hep-th]].
  
\bibitem{Goldfeld} D. Goldfeld, \emph{ A Simple Proof of Siegel's Theorem}, Proc. Nat. Acad. Sci. USA, {\bf 71} No. 5 (1974), 1055.

\bibitem{Goldfeld2} D. Goldfeld, \emph{Gauss' class number problem for imaginary quadratic fields},
Bull. Amer. Math. Soc.  {\bf 13}, (1985),  23-37.

\bibitem{Goldstein} 
  K.~Goldstein, N.~Iizuka, R.~P.~Jena and S.~P.~Trivedi,
  \emph{Non-supersymmetric attractors},
  Phys.\ Rev.\ D {\bf 72} (2005), 124021 [hep-th/0507096].


\bibitem{Gomes} J. Gomes, \emph{Exact holography and black hole entropy in $N=8$ and $N=4$ string theory}, JHEP {\bf 07} (2017), 022 [arXiv:1511.07061 [hep-th]].

\bibitem{Heegner}
K. Heegner, \emph{Diophantische Analysis und Modulfunktionen}, Mathematische Zeitschrift, {\bf 56} (3) (1952), 227--253.


\bibitem{Heilbronn} H. Heilbronn, \emph{On the class-number in imaginary quadratic fields}, Quart. J. Math. Oxford Ser. {\bf 5} (1934), 150--160.

\bibitem{Hohn1} G. H\"ohn, ``Conformal Designs based on Vertex Operator Algebras,'' arXiv:math/0701626.

\bibitem{Hohn2} G. H\"ohn, ``Selbstduale Vertexoperatorsuperalgebren und das Babymonster (Self-dual Vertex Operator Super Algebras and the Baby Monster),'' Ph.D. thesis  (Bonn 1995), Bonner Mathematische Schriften 286 (1996), 1-85, arXiv:0706.0236.

\bibitem{KachruTripathy} S. Kachru and A. Tripathy, \emph{Black holes and Hurwitz class numbers},
Int. J. Modern Physics D{\bf{26}}, 12 (2017),  1742003 [arXiv:1705.06295 [hep-th]].

\bibitem{Kawai} T. Kawai, Y. Yamada and S.K. Yang, \emph{Elliptic genera and N=2 superconformal field theory}, Nucl. Phys. B{\bf 414} (1994), 191--212 [hep-th/9306096].

\bibitem{KohnenOno} W. Kohnen and K. Ono, 
\emph{Indivisibility of class numbers of imaginary quadratic fields and orders of Tate-Shafarevich groups of elliptic curves with complex multiplication}, Inventiones Mathematicae {\bf 135} (1999),
387--398.

\bibitem{Juan} 
  J.~M.~Maldacena,
  \emph{The Large N limit of superconformal field theories and supergravity},
  Int.\ J.\ Theor.\ Phys.\  {\bf 38}, 1113 (1999)
  [Adv.\ Theor.\ Math.\ Phys.\  {\bf 2}, 231 (1998)]
  [hep-th/9711200].
  
  \bibitem{Manschot} 
  J.~Manschot, \emph{On the space of elliptic genera},
  Commun.\ Num.\ Theor.\ Phys.\  {\bf 2}, 803 (2008)
  doi:10.4310/CNTP.2008.v2.n4.a4
  [arXiv:0805.4333 [hep-th]].
  

\bibitem{Moore3} G. W. Moore, \emph{Attractors and arithmetic},
arXiv:hep-th/9807056.

\bibitem{Moore1} G. W. Moore, \emph{Arithmetic and attractors}, arXiv:hep-th/9807087.

\bibitem{Moore2} G. W. Moore, \emph{Strings and arithmetic}, arXiv:hep-th/0401049.

\bibitem{WittenNote} K. Ono and L. Rolen, \emph{On Witten's extremal partition functions}, submitted.

\bibitem{Polak} J. Polak, unpublished.

\bibitem{Schoof} R. Schoof, \emph{Nonsingular plane cubic curves over finite fields},
Journal of Combinatorial Theory, Series A
{\bf 46}, Issue 2, (1987), 183--211.

\bibitem{Sengen} A. Sen, \emph{Black hole solutions in heterotic
string theory on a torus}, Nucl. Phys. B {\bf 440} (1995), 421--440 [hep-th/9411187].

\bibitem{Senlect} 
  A.~Sen,
  \emph{Black Hole Entropy Function, Attractors and Precision Counting of Microstates},
  Gen.\ Rel.\ Grav.\  {\bf 40} (2008), 2249
  [arXiv:0708.1270 [hep-th]].


\bibitem{ShiodaInose} T. Shioda and H. Inose, \emph{On singular $K3$ surfaces}, in Complex analysis and algebraic
geometry, Cambridge University Press, Cambridge, 1977.

\bibitem{Siegel} C. L. Siegel, \emph{\"Uber die Classenzahl quadratischer Zahlk\"orper}, Acta Arithmetica {\bf 1} (1935), 83--86.

\bibitem{Sound} K. Soundararajan, \emph{Divisibility of Class Numbers of Imaginary Quadratic Fields}, Journal of the London Mathematical Society, {\bf 61}, Issue 3 (2000), 681--690.

\bibitem{Stark} H. M. Stark, \emph{A Complete Determination of the Complex Quadratic Fields of Class Number One}, Michigan Math. J. {\bf 14} (1967), 1--27.
  
\bibitem{Strominger} 
  A.~Strominger and C.~Vafa,
  \emph{Microscopic origin of the Bekenstein-Hawking entropy},
  Phys.\ Lett.\ B {\bf 379} (1996), 99
  [hep-th/9601029].

\bibitem{Tripathy}
  P.~K.~Tripathy and S.~P.~Trivedi,
  \emph{Non-supersymmetric attractors in string theory},
  JHEP {\bf 0603} (2006), 022 
  [hep-th/0511117].

\bibitem{Vafa}
  C.~Vafa,
  \emph{Instantons on D-branes},
  Nucl.\ Phys.\ B {\bf 463} (1996), 435
  [hep-th/9512078].


\bibitem{Watkins} M. Watkins, \emph{ Class numbers of imaginary quadratic fields, Mathematics of Computation}, {\bf 73} (2004), 907--938.

\bibitem{Witten} 
  E.~Witten,
  \emph{Three-Dimensional Gravity Revisited},
  arXiv:0706.3359 [hep-th].

\bibitem{Zagier} D. Zagier, \emph{Traces des op\'erateurs de Hecke}, 
S\'eminaire Delange-Pisot-Poitou 1975-1976, Expos\'e No. 23, 12 pages,
reprinted (in English translation) as: {\it The Eichler-Selberg trace formula on $\mathrm{SL}_2(\Z)$}, 
Appendix to S. Lang, {\it Introduction to Modular Forms}, Grundlehren d. math. Wiss. 222, Springer-Verlag, Berlin-Heidelberg-New York (1976) 44-54.


\end{thebibliography}
\end{document}